\documentclass[microtype]{jloganal}

\author[H Towsner]{Henry Towsner}
\givenname{Henry}
\surname{Towsner}
\address{Department of Mathematics\\
University of Connecticut}
\email{htowsner@gmail.com}
\urladdr{http://www.math.uconn.edu/$\sim$towsner}
\title{A Correspondence Principle for the Gowers Norms}
\date{\today}

\arxivreference{math.DS/0905.0493}
\arxivpassword{8ziuu}

\newtheorem{lemma}{Lemma}[section]
\newtheorem{theorem}{Theorem}[section]
\newtheorem{definition}{Definition}[section]

\volumenumber{}
 \issuenumber{}
\publicationyear{}
\papernumber{} 
\startpage{} 
\endpage{}
\doi{} 
\MR{} 
\Zbl{} 
\received{} 
\revised{} 
\accepted{} 
\published{}
\publishedonline{} 
\proposed{} 
\seconded{} 
\corresponding{} 
\editor{} 
\version{}

\begin{document}

\begin{abstract}
  The Furstenberg Correspondence shows that certain ``local behavior'' of dynamical system is equivalent to the behavior of sufficiently large finite systems.  The Gowers uniformity norms, however, are not local in the relevant sense.  We give a modified correspondence in which the Gowers norm is preserved.  This extends to the integers a similar result by Tao and Zielger on finite fields.
\end{abstract}

\maketitle

\section{Introduction}
Informally speaking, the Furstenberg Correspondence \cite{furstenberg77,furstenberg82} shows that the ``local behavior'' of a dynamical system is controlled by the behavior of sufficiently large finite systems.  By the local behavior of a dynamical system $(X,\mathcal{B},\mu,G)$, we mean the properties which can be stated using finitely many actions of $G$ and the integral given by $\mu$\footnote{A more precise version of this notion would be to say that the local behavior consists of the $\Pi_2$ formulas in an appropriate extension of the language of arithmetic.}.  By a finite system, we just mean $(S,\mathcal{P}(S),c,G)$ where $G$ is a infinite group, $S$ is a finite quotient of $G$, and $c$ is the counting measure $c(A):=\frac{|A|}{|S|}$.

The most well known example of such a property is the ergodic form of Szemer\'edi's Theorem:
\begin{quote}
  For every $k$, every $\epsilon>0$, and every $L^\infty$ function $f$, if $\int f d\mu\geq\epsilon$ then there is some $n$ such that $\int\prod_{j=0}^{k-1}T^{-jn}f d\mu>0$.
\end{quote}
The Furstenberg Correspondence shows that this is equivalent to the following statement of Szemer\'edi's Theorem:
\begin{quote}
  For every $k$ and every $\epsilon>0$, there is an $N$ and a $\delta>0$ such that if $m\geq N$ and $f:[0,m-1]\rightarrow[-1,1]$ is such that $\int f dc\geq\epsilon$ then there is some $n$ such that $\int \prod_{j=0}^{k-1}f(x+jn) dc(x)\geq\delta$.
\end{quote}
In general, the Furstenberg Correspondence states that, given a sequence of functions on increasingly large finite systems, a single function on a single infinite system can be given with the property that suitable calculations are controlled by the limit of the value of analogous calculations in the finite systems.

Recent work by Austin and Tao \cite{austin08}, Elek and Szegedy \cite{elek07}, Elek \cite{elek08}, Tao \cite{tao07}, and the author \cite{towsner:MR2529651,towsner08} has extended this correspondence both to other specific properties and to more general formulations.  These methods are not adequate, however, for the study of the uniformity norms introduced for finite systems by Gowers in \cite{gowers01} and for infinite systems by Host and Kra in \cite{host05}.  While there are strong reasons for believing that functions on finite systems with small Gowers norm should correspond to functions on infinite systems with small Gowers-Host-Kra norm, these norms are not local.
In particular, the ordinary correspondence may place a sequence of highly $k$-uniform functions (that is, functions with $||\cdot||_{U^k}$ going to zero) in correspondence with a function with large $U^k$ norm.


In \cite{tao08}, Tao and Ziegler give a variant of the correspondence principle which preserves the $U^k$ norms when the group $G$ is vector space over a finite field, $\mathbb{F}^\omega_p$.  Their argument, however, takes advantage of group theoretic properties of $\mathbb{F}_p$, and does not immediately extend to other groups.

In this paper we give a similar correspondence for arbitrary countable Abelian groups.  While there is no theoretical obstacle to giving the construction explicitly in a style similar to \cite{tao08}, the resulting argument would be quite unwieldy.  Roughly speaking, where Tao and Ziegler can choose representative transformations randomly and expect that almost all choices suffice, here we have to choose particular transformations.  It is much more convenient to do the work of choosing the correct transformations in an infinitary ergodic setting; the price is that we use an argument from nonstandard analysis to give a highly infinitary system acted on by a very large group, and then use ergodic methods to reduce the group down to something more manageable.  For the sake of readers unfamiliar with nonstandard analysis, we isolate its use in a single lemma.

In \fullref{section:subgroup} we lay the ergodic-theoretic groundwork for the correspondence, and in \fullref{section:nonstandard} we give the correspondence argument itself.

The author is grateful to the Mathematical Sciences Research Institute's semester program in Ergodic Theory and Additive Combinatorics, during which the ideas in this paper were developed.

\section{Choosing a Good Subgroup}
\label{section:subgroup}
Because of the nature of the intermediate object which will be produced by the nonstandard argument in \fullref{section:nonstandard}, we want to work with a fairly general notion of a dynamical system.

\begin{definition}
  A \emph{dynamical system} consists of a probability measure space $(X,\mathcal{B},\mu)$ together with an Abelian group $G$, an action of $G$ on $X$ such that for each $g\in G$, the action $T_g:X\rightarrow X$ is measurable, and a finitely additive $G$-invariant probability measure space $(G,\mathcal{C},\lambda)$.
\end{definition}
We do not require that the action of $G$ on $X$ be a measurable function from $G\times X$ to $X$, since we need actions where this is not true.  Instead we only ask for the weaker condition that Fubini's Theorem holds.

The most common case is where $G$ is countable and $\mathcal{C}$ is the powerset of $G$ (which is possible since $\lambda$ is only required to be finitely additive).  
\begin{definition}
By a \emph{discrete group}, we mean a countable group $G$ together with its power set, viewed as a measure space.
\end{definition}
This case is common enough that when $G$ is countable, we will write (as we did in the introduction) $(X,\mathcal{B},\mu,G)$ to mean the dynamical system $(X,\mathcal{B},\mu),(G,\mathcal{P}(G),\lambda),\{T_g\}_{g\in G}$ where $\lambda$ and $\{T_g\}$ are implicit or given by the context.

\begin{definition}
  Given a bounded function $f$ and a group $G$, define $\mathcal{F}(f,G)$ to be the collection of functions containing $f$ and the function constantly equal to $1$, and closed under pairwise sums, pairwise multiplication, scalar multiplication by a rational, and shifts from $G$. 
\end{definition}
 Clearly $\mathcal{F}(f,G)$ is countable as long as $G$ is.

\begin{definition}
When $(X,\mathcal{B},\mu)$ is a probability measure space and $f:X\rightarrow\mathbb{R}$ is bounded and measurable, we say $f$ is \emph{weakly Fubini} for $(G,\mathcal{C},\lambda)$ if for every $x\in X$, the function $g\mapsto f(T_gx)$ is measurable with respect to $\mathcal{C}$ and $x\mapsto \int f(T_gx) d\lambda$ is measurable with respect to $\mathcal{B}$.

  We say $f$ is \emph{Fubini} if 
every function in $\mathcal{F}(f,G)$ is weakly Fubini.
\end{definition}
The requirement that the condition hold for every $x$ could be weakened to almost every $x$ without much trouble.

In this context, the Mean Ergodic Theorem can be taken to be the following:
\begin{lemma}
  If $f$ is Fubini for $(G,\mathcal{C},\lambda)$ and $\mathcal{I}(G)$ is the collection of sets invariant under $T_g$ for every $g\in G$ then 
\[\int \left[E(f\mid\mathcal{I}(G))(x)-\int f(T_gx) d\lambda\right]^2 d\mu=0.\]
\end{lemma}

\begin{lemma}
If $G$ is countable, every measurable function $f$ on $X$ is weakly Fubini for $(G,\mathcal{P}(G),\lambda)$.
\end{lemma}
\begin{proof}
Measurability of $g\mapsto f(T_gx)$ follows because every subset of $G$ is measurable.  The standard construction of a F\o{}lner sequence \cite{namioka:MR0180832} gives a sequence of subsets $I_N$ of $G$ so that $\int f(T_gx)d\lambda=\lim_{N\rightarrow\infty, N\in U}\frac{1}{|I_N|}\sum_{g\in I_N}f(T_gx)$ for some infinite set $U$.  Since each of the functions $x\mapsto \frac{1}{|I_N|}\sum_{g\in I_N}f(T_gx)$ is measurable, so is the limit. 
\end{proof}

Similarly, we may extend a group by a single element (or more precisely, by the discrete group generated by that element) while preserving the Fubini property.
\begin{definition}
  Let $H$ be a subgroup of $G$ and suppose $f$ is Fubini for $(G,\mathcal{C},\lambda)$.  For $g\in G$, define $H'_g$ to be the subgroup of $G$ generated by $H\cup\{g\}$.  Taking $\pi:H\times\mathbb{Z}\rightarrow H'_g$ to be the homomorphism given by $\pi(h,n):=h\cdot g^n$, any finitely-additive $H$-invariant measure $(H,\mathcal{D},\nu)$ may be extended to a measure on $H'_g$ by taking the $\pi$-image of $(H\times \mathbb{Z},\mathcal{D}\times\mathcal{P}(\mathbb{Z}),\nu\times\sigma)$ where $\sigma$ is an arbitrary finitely-additive $\mathbb{Z}$-invariant measure.
\end{definition}
Note that the choice of measure on $H'_g$ is not canonical, but the Mean Ergodic Theorem tells us that the choice will not matter.

\begin{lemma}
	Suppose $H\subseteq G$, $f$ is bounded and Fubini for $(G,\mathcal{C},\lambda)$ and $(H,\mathcal{D},\nu)$, and the inclusion map of $H$ in $G$ is measurable as a function from $(H,\mathcal{D})$ to $(G,\mathcal{C})$.  Then for any $g\in G$,  $f$ is Fubini for $H'_g$.
\end{lemma}
\begin{proof}
For any $f'\in\mathcal{F}(f,H)$, any measurable set $S$ and any $n$, $\{h\mid f'(T_h(T_{g^n}x))\in S\}=\{h\mid f'(T_{h\cdot g^n}x)\in S\}$ is measurable (in $\mathcal{D}$).  Therefore $\{(h,n)\mid f'(T_{h\cdot g^n}x)\in S\}=\bigcup_n\{h\mid f'(T_{h\cdot g^n}x\in S)\}\times\{n\}$ is measurable as well.

By the same argument as above, there is a set $U$ such that
\[\int f'(T_hx)d(\nu\times\sigma)(h)=\lim_{N\rightarrow\infty,N\in U}\frac{1}{|I_N|}\sum_{i\in I_n}\int f'(T_hT_{g^n}x)d\nu(h)\]
and since $f'\circ T_{g^n}$ is weakly Fubini, each $\frac{1}{|I_N|}\sum_{n\in I_n} \int f'(T_h T_{g^n} x) d\nu(h)$ is measurable, and therefore the limit is as well.
\end{proof}

The following definition and the basic properties of such norms are taken from \cite{host05}.
\begin{definition}[Gowers-Host-Kra Norms]
  Define $X^{[k]}:=X^{2^k}$, $\mathcal{B}^{[k]}:=\mathcal{B}^{2^k}$, and for any transformation $T$ on $(X,\mathcal{B},\mu)$, define $T^{[k]}:=\bigotimes_{\omega\in \{0,1\}^k}T$.

  If $G$ is an Abelian group acting on $(X,\mathcal{B},\mu)$, define $\mu^{[0]}(G):=\mu$, $\mathcal{I}^{[k]}(G)$ to be the collection of sets in $\mathcal{B}^{[k]}$ invariant under $T^{[k]}_g$ for each $g\in G$, and $\mu^{[k+1]}(G)$ to be the relative joining of $\mu^{[k]}$ with itself over $\mathcal{I}^{[k]}(G)$.

For any $L^\infty$ function $f:X\rightarrow\mathbb{R}$, define
\[||f||_{U^k(G)}:=\left(\int \bigotimes_{\omega\in\{0,1\}^k} f d\mu^{[k]}(G)\right)^{1/2^k}.\]
\end{definition}
Note that this definition depends on $G$, and the action of $G$ on $X$, but not on a particular measure space on $G$.  Despite the name, $||\cdot||_{U^k(G)}$ is generally only a semi-norm.  (There is a similar norm for complex valued functions, the only difference being replacing $f$ by the complex conjugate of $f$ in some cases; the arguments in this paper go through unchanged for complex valued functions, so we only discuss the notationally simpler real valued version.)

\begin{lemma}
  If $f$ is Fubini for $(G,\mathcal{C},\lambda)$ then
\[||f||^{2^k}_{U^k(G)}=\iint \prod_{\omega\in\{0,1\}^k}T_{\omega\cdot\vec g}f d\lambda^k(\vec g)d\mu.\]
\end{lemma}
\begin{proof}
  We show by induction on $k$ that
\[\int \bigotimes_{\omega\in\{0,1\}^k} f_\omega d\mu^{[k]}=\iint \prod_{\omega\in\{0,1\}^k}T_{\omega\cdot\vec g}f_\omega d\mu d\lambda^k(\vec g).\]
For $k=0$, this is immediate.  Assume the claim holds for $k$.  Then
\begin{align*}
\int \bigotimes_{\omega\in\{0,1\}^{k+1}} f_\omega d\mu^{[k+1]}
=&\int E(\bigotimes_{\omega\in\{0,1\}^k} f_{0^\frown\langle \omega\rangle}\mid\mathcal{I}^{[k]}(G))\cdot\\
&\ \ \ \ \ \ \ \ \ \ \ \ \ \ \ \ E(\bigotimes_{\omega\in\{0,1\}^k} f_{1^\frown\langle \omega\rangle}\mid\mathcal{I}^{[k]}(G)) d\mu^{[k]}(G)\\
=&\iint \bigotimes_{\omega\in\{0,1\}^k} f_{0^\frown\langle\omega\rangle} T_g f_{1^\frown\langle \omega\rangle} d\lambda d\mu^{[k]}(G)\\
=&\iint \bigotimes_{\omega\in\{0,1\}^k}f_{0^\frown\langle\omega\rangle} T_g f_{1^\frown\langle \omega\rangle} d\mu^{[k]}(G)d\lambda\\
=&\iiint \prod_{\omega\in\{0,1\}^k}T_{\omega\cdot\vec g}(f_{0^\frown\langle\omega\rangle} T_g f_{1^\frown\langle \omega\rangle}) d\mu d\lambda^k(\vec g) d\lambda(g)\\
=&\iint \prod_{\omega\in\{0,1\}^{k+1}}T_{\omega\cdot\vec g} f_\omega d\mu d\lambda^{k+1}(\vec g)\\
\end{align*}
\end{proof}

The following property is easily seen by induction:
\begin{lemma}
  If $H$ is a subgroup of $G$ then $||f||_{U^k(G)}\leq||f||_{U^k(H)}$.
\end{lemma}

It will be convenient to use a slight generalization of the $U^k$ norm, in which a different group is used at the top-most level.
\begin{definition}
  Let $H,G$ be groups.  Then $\mu^{[0]}(G,H):=\mu$, $\mathcal{I}^{[k]}(G,H)$ is the space of sets $B\in\mathcal{B}^{[k]}$ such that $\mu^{[k]}(G)(B\bigtriangleup (T_h^{[k]})^{-1}B)=0$ for each $h\in H$, and $\mu^{[k+1]}(G,H)$ is the relative joining of $\mu^{[k]}(G)$ with itself over $\mathcal{I}^{[k]}(G,H)$.

Similarly, $||f||_{U^k(G,H)}:=\left(\int \bigotimes_{\omega\in\{0,1\}^k} f d\mu^{[k]}(G,H)\right)^{1/2^k}$.
\end{definition}

\begin{lemma}
  If $H$ and $H'$ are subgroups of $G$ and $H$ is a subgroup of $H'$ then $||f||_{U^k(G,H')}\leq||f||_{U^k(G,H)}$.
\end{lemma}

\begin{theorem}
  Let $(X,\mathcal{B},\mu)$, $(G,\mathcal{C},\lambda)$ be a dynamical system, let $H$ be a subgroup of $G$, and let $(H,\mathcal{C}',\lambda')$ be given so the inclusion of $H$ in $G$ is measurable.  Let $f$ be everywhere bounded by $1$, let $f$ be Fubini for both $G$ and $H$, and suppose that $||f||^{2^{k+1}}_{U^{k+1}(G,H)}=||f||_{U^{k+1}(G)}^{2^{k+1}}+\epsilon$ with $\epsilon>0$.  Then there is a $g\in G$ such that, $||f||^{2^{k+1}}_{U^{k+1}(G,H'_g)}\leq ||f||^{2^{k+1}}_{U^{k+1}(G)}+3\epsilon/4$.
\end{theorem}
\begin{proof}
  Note that a similar claim for the $U^0$ norm would be trivial, since the premise could never hold (the $U^0$ norm is independent of $G$).  Observe that
\[||f||^{2^{k+1}}_{U^{k+1}(G)}=\int E(\bigotimes_{\omega\in\{0,1\}^k} f\mid\mathcal{I}^{[k]}(G))^2 d\mu^{[k]}(G).\]
Setting $f':=E(\bigotimes f\mid\mathcal{I}^{[k]}(G,H))$, this quantity is equal to
\[\int E(f'\mid\mathcal{I}^{[k]}(G))^2 d\mu^{[k]}(G)=||E(f'\mid\mathcal{I}^{[k]}(G))||_{L^2(\mu^{[k]}(G))}^2.\]
Suppose that for every $g\in G$, $||f'-\lim_{n\rightarrow\infty}\frac{1}{n}\sum_{i\leq n}(T_g^{[k]})^if'||_{L^2(\mu^{[k]}(G))}<\sqrt{\epsilon}/2$.  Then also $||f'-E(f'\mid\mathcal{I}^{[k]}(G))||^2_{L^2(\mu^{[k]}(\Gamma))}\leq\epsilon/2$, which implies that
\[\left| ||f||_{U^{k+1}(G)}^{2^{k+1}}-||f||^{2^{k+1}}_{U^{k+1}(G,H)}\right|=\left| ||E(f'\mid\mathcal{I}^{[k]}(G))||_{L^2}-||f'||_{L^2}\right|\leq\epsilon/2,\]
contradicting the assumption.

So choose $g$ such that $||f'-\lim_{n\rightarrow\infty}\frac{1}{n}\sum_{i\leq n}(T_g^{[k]})^i||\geq\sqrt{\epsilon}/2$.  It follows that $||f'||^2-||E(f'\mid\mathcal{I}^{[k]}(H'_g)||^2=||f'-E(f'\mid\mathcal{I}^{[k]}(H'_g))||^2\geq \epsilon/4$, and therefore that
\[||f||^{2^{k+1}}_{U^{k+1}(G,H'_g)}=||E(f'\mid H'_g)||^2\leq ||f||^{2^{k+1}}_{U^{k+1}(G)}+\frac{3}{4}\epsilon.\]
\end{proof}

\begin{lemma}
  Let $(X,\mathcal{B},\mu)$, $(G,\mathcal{C},\lambda)$ be a dynamical system and, let $\Gamma$ be a discrete subgroup of $G$, let $f$ be bounded and Fubini for $(G,\mathcal{C},\lambda)$.  Then there is a discrete subgroup $H$ of $G$ containing $\Gamma$ such that $||F||_{U^k(G,H)}=||F||_{U^k(G)}$ for every $F\in\mathcal{F}(f,H)$.
\end{lemma}
\begin{proof}
  We will construct $H$ so that there is a natural map $\pi:\mathbb{Z}^\omega\rightarrow H$ (by $\mathbb{Z}^\omega$ we mean the product of countably many copies of $\mathbb{Z}$, itself a countable set).  Then we may choose a sequence of pairs $(\epsilon_n,F_n)$ where $\epsilon_n$ is a rational in $(0,1)$ and $F_n$ is a code for an element of $\mathcal{F}(f,H)$ so that each such pair appears at some point in this sequence.  We set $H_0:=\Gamma$, and for each $n$, set $H_{n+1}:=(H_n)'_g$ where $g$ is chosen so that $||F_n||_{U^k(G,H_{n+1})}\leq||F_n||_{U^k(G)}+\epsilon_n$.  Take $H:=\bigcup_{n<\omega}H_n$.  Then for every $F\in\mathcal{F}(f,H)$ and every $\epsilon>0$, for sufficiently large $n$, $||F||_{U^k(G,H_n)}\leq||F||_{U^k(G)}+\epsilon$.  $H$ is a subgroup of $G$, so $||F||_{U^k(G)}\leq||F||_{U^k(G,H)}$, and therefore $||F||_{U^k(G)}=||F||_{U^k(G,H)}$.
\end{proof}

\begin{theorem}
  Let $(X,\mathcal{B},\mu)$, $(G,\mathcal{C},\lambda)$ be a dynamical system, let $\Gamma$ be a discrete subgroup of $G$, and let $f$ be bounded and Fubini for $(G,\mathcal{C},\lambda)$.  Then there is a discrete subgroup $H$ of $G$ containing $\Gamma$ such that $||F||_{U^k(H)}=||F||_{U^k(G)}$ for every $F\in\mathcal{F}(f,H)$.
\end{theorem}
\begin{proof}
Let $H$ be as given in the preceeding lemma and proceed by induction on $k$.  For $k=0$, this is trivial.  Assume the result holds for $k$.  Then
\[||F||_{U^{k+1}(H)}^{2^{k+1}}=\int E(\bigotimes F\mid\mathcal{I}^{[k]}(\mathbb{Z}^\omega))^2 d\mu^{[k]}(H).\]
For any $\epsilon>0$, we may choose $i$ so that this is within $\epsilon$ of
\[\int E(\bigotimes F\mid\mathcal{I}^{[k]}(H_i))^2 d\mu^{[k]}(H)=\int \bigotimes FT_hF d\mu^{[k]}(H).\]
For every $h\in H_i$, 
\[\int \bigotimes F\cdot T_hF d\mu^{[k]}(H)=||FT_hF||^{2^k}_{U^k(H)}\]
and since $FT_hF\in\mathcal{F}(f,H)$, by IH $||FT_hF||^{2^k}_{U^k(H)}=||FT_hF||^{2^k}_{U^k(G)}$.
It follows that
\[\left|||F||^{2^{k+1}}_{U^{k+1}(H,H_i)}-||F||^{2^{k+1}}_{U^{k+1}(G,H_i)}\right|<\epsilon.\]
But for sufficiently large $i$, $||F||^{2^{k+1}}_{U^{k+1}(G,H_i)}$ is arbitrarily close to $||F||^{2^{k+1}}_{U^{k+1}(G)}$.  So, taking the limit as $i\rightarrow\infty$, we have
\[||F||^{2^{k+1}}_{U^{k+1}(H)}=||F||^{2^{k+1}}_{U^{k+1}(G)}.\]
\end{proof}

\section{A Correspondence Principle}
\label{section:nonstandard}

To set up the appropriate analogy between different dynamical systems, we need the notion of a representative of an element of $\mathcal{F}(\cdot, G)$.
\begin{definition}
  Let $\mathcal{F}(G)$ be a set of symbols defined inductively by:
  \begin{itemize}
  \item $\mathfrak{c}\in\mathcal{F}(G)$
  \item $\mathfrak{1}\in\mathcal{F}(G)$
  \item If $\mathfrak{f},\mathfrak{g}\in\mathcal{F}(G)$ then $\mathfrak{f}+\mathfrak{g}$ and $\mathfrak{f}\cdot \mathfrak{g}$ belong to $\mathcal{F}(G)$
  \item If $h\in G$ and $\mathfrak{f}\in\mathcal{F}(G)$ then $\mathfrak{T}_h\mathfrak{f}\in\mathcal{F}(G)$
  \item If $\mathfrak{f}\in\mathcal{F}(G)$ and $q$ is a rational then $q\mathfrak{f}\in\mathcal{F}(G)$
  \end{itemize}

  If $\mathfrak{f}\in\mathcal{F}(G)$ and $f$ is a bounded Fubini function in a dynamical system, we define $\mathfrak{f}(f)$ recursively by:
  \begin{itemize}
  \item $\mathfrak{c}(f):=f$
  \item $\mathfrak{1}(f):=1$ (the function constantly equal to $1$)
  \item $(\mathfrak{f}+\mathfrak{g})(f):=\mathfrak{f}(f)+\mathfrak{g}(f)$
  \item $(\mathfrak{f}\cdot\mathfrak{g})(f):=\mathfrak{f}(f)\cdot\mathfrak{g}(f)$
  \item $(\mathfrak{T}_h\mathfrak{f})(f):=T_h(\mathfrak{f}(f))$
  \item $(q\mathfrak{f})(f):=q\cdot(\mathfrak{f}(f))$
  \end{itemize}
\end{definition}
It is easy to see that $g\in\mathcal{F}(f,G)$ iff there is a $\mathfrak{f}\in\mathcal{F}(G)$ such that $g=\mathfrak{f}(f)$.

\begin{lemma}
  Let $\Gamma$ be a countable Abelian group, let $\mathcal{N}$ be an infinite set of integers, and for each $N\in\mathcal{N}$, let $S_N$ be a finite quotient of $\Gamma$, $\pi_N:\Gamma\rightarrow S_N$, with $|S_N|\rightarrow\infty$.  Let $f_N:S_N\rightarrow[-1,1]$ be given.  There is a dynamical system $(X,\mathcal{B},\mu)$, $(G,\mathcal{C},\lambda)$ with $\lambda$ $\sigma$-additive, a homomorphism $\pi:\Gamma\rightarrow G$, and a measurable Fubini function $f:X\rightarrow[-1,1]$ such that for any $\mathfrak{f}\in\mathcal{F}(\Gamma)$,
\[\liminf_{N\in\mathcal{N}}||\mathfrak{f}(f_N)||_{U^k(S_N)}\leq||\mathfrak{f}(f)||_{U^k(G)}\leq\limsup_{N\in\mathcal{N}}||\mathfrak{f}(f_N)||_{U^k(S_N)}\]
for each $k$.

Additionally, if for every $g,h\in \Gamma$, $g\neq h$, $\pi_N(g)\neq\pi_N(h)$ except for finitely many $N\in\mathcal{N}$ then $\pi(g)\neq\pi(h)$.
\end{lemma}
\begin{proof}
  Fix a non-principal ultrafilter $U$ and form the nonstandard extension of a universe containing the sequences $\mathcal{N}$, $\langle f_N\rangle$, $\langle S_N\rangle$.  The sequence $\mathcal{N}$ codes a nonstandard integer $a$ and $G=S_a$ is a hyperfinite Abelian group.  By the Loeb measure construction, the internal subsets of $G$ may be extended to a $\sigma$-algebra $\mathcal{C}$ on $G$, the internal counting measure on $G$ may be extended to a $\sigma$-additive measure $\lambda$, and the sequence $\langle f_N\rangle$ represents an internal function $F:\Gamma\rightarrow[-1,1]^*$.  The function $f=st\circ F$ is then a measurable function from $G$ to $[-1,1]$.

This same measure space is also $(G,\mathcal{C},\lambda)$, with $G$ acting on itself by the group action.  Since $F$ is internal, for any $x\in S_a$, the functions $g\mapsto f(T_gx)$ and $x\mapsto \int f(T_gx)d\lambda$ are the result of applying the standard part operation to internal functions, and are therefore measurable, and the same applies to any element of $\mathcal{F}(f,G)$.  So $f$ is Fubini.  The embedding $\pi:\Gamma\rightarrow G$ is simply the embedding represented by the sequence $\langle \pi_N\rangle$.

The final clause follows from transfer.  For instance, if $\liminf_{N\in\mathcal{N}}||f_N||_{U^k(S_N)}\geq\alpha$ then for each $\epsilon>0$ and all but finitely many $N$ in $\mathcal{N}$,
\[\frac{1}{|S_N|^{k+1}}\sum_{x\in S_N}\sum_{\vec g\in S_N^k}\prod_{\omega\in\{0,1\}^k}f_N(x+\vec g\cdot\omega)>\alpha-\epsilon\]
and therefore
\[st\left(\frac{1}{|S_a|^{k+1}}\sum_{x\in S_a}\sum_{\vec g\in S_a^k}\prod_{\omega\in\{0,1\}^k}F(x+\vec g\cdot\omega)\right)\geq\alpha-\epsilon\]
and therefore
\[\iint \prod_{\omega\in\{0,1\}^k}T_{\vec g\cdot\omega}f d\mu d\lambda^k\geq\alpha-\epsilon.\]
Applying this argument for arbitrary $\mathfrak{f}\in\mathcal{F}(\Gamma)$, and the analogous argument for the upper bound, gives the claim.
\end{proof}

As shown in the previous section, there is a discrete subgroup $H$ of $G$ containing $\Gamma$ such that for each $F\in\mathcal{F}(f,H)$, $||F||_{U^k(H)}=||F||_{U^k(G)}$.  Putting this together, we obtain the following theorem:
\begin{theorem}
  Let $\Gamma$ be a countable Abelian group, let $\mathcal{N}$ be an infinite set of integers, and for each $N\in\mathcal{N}$, let $S_N$ be a finite quotient of $\Gamma$, $\pi_N:\Gamma\rightarrow S_N$, with $|S_N|\rightarrow\infty$.  Let $f_N:S_N\rightarrow[-1,1]$ be given.  Then there is a dynamical system $(X,\mathcal{B},\mu,H)$ and an $L^\infty$ function $f:X\rightarrow\mathbb{C}$ such that for any $\mathfrak{f}\in\mathcal{F}(\Gamma)$, 
\[\liminf_{N\in\mathcal{N}}||\mathfrak{f}(f_N)||_{U^k(S_N)}\leq||\mathfrak{f}(f)||_{U^k(H)}\leq\limsup_{N\in\mathcal{N}}||\mathfrak{f}(f_N)||_{U^k(S_N)}\]
for each $k$.
\end{theorem}



\end{document}